\documentclass[12pt]{article}
\usepackage{ifpdf}
\usepackage {graphics}
\usepackage{amsmath,mathptmx,amssymb,bm}
\newtheorem{theorem}{Theorem}

\newtheorem{definition}{Definition}

\begin{document}


\title{Equivalence transformations and conservation laws for a generalized variable-coefficient Gardner equation}


\author{R. de la Rosa${}^{a}$, M.L. Gandarias${}^{b}$, M.S. Bruz\'on${}^{c}$\\
 ${}^{a}$ Universidad de C\'adiz, Spain  (e-mail: rafael.delarosa@uca.es). \\
 ${}^{b}$ Universidad de C\'adiz, Spain (e-mail:  marialuz.gandarias@uca.es). \\
 ${}^{c}$ Universidad de C\'adiz, Spain  (e-mail:  m.bruzon@uca.es). \\  
}

\date{}
 
\maketitle

\begin{abstract}

In this paper we study the generalized variable-coefficient Gardner equations of the form $u_t + A(t)u^n\,u_x+ C(t)\,u^{2n}u_x + B(t)\,u_{xxx} + Q(t)\,u =0$. This class broadens out many other equations previously considered: Johnpillai and Khalique (2010), Molati and Ramollo (2012) and Vaneeva, Kuriksha and Sophocleous (2015). Equivalence group of the class under consideration has been constructed which permit an exhaustive study and a simple and clear formulation of the results. Some conservation laws are derived for the nonlinearly self-adjoint equations, based on differential substitutions, and by using the direct method of the multipliers.\\

 
\end{abstract}



\maketitle

\section{Introduction}
\label{}

Recent developments in the field of partial differential equations (PDEs) have led researchers to focus its efforts on the study of PDEs with variable coefficients, particularly in nonlinear equations with variable coefficients. These equations describe many nonlinear phenomena more realistically than their constant coefficients counterparts. However, the study of variable coefficient equations seems rather difficult. Nowadays, the research on variable coefficient equations primarily encompasses the study of the symmetries of PDEs, the determination of integrability conditions, the construction of conservation laws or the obtaining of exact solutions.\\

Nonlinear evolution equations play an important role in the field of nonlinear dynamics. Among them, we emphasise the KdV equation as well as their generalizations. KdV equation governs different physical processes, for instance, the dynamics and the physics of shallow waters. The problem lies in the fact that KdV equation is a quite simple model to analyse these phenomena more precisely. Thus, it must be considered generalizations of KdV equation which involve more than one nonlinear term, the Gardner equation is an example. The Gardner equation is also known as combined KdV-mKdV equation, this is a useful model for the description of wave phenomena in plasma and solid state and internal solitary waves in shallow waters.\\

In recent years, several works have been dedicated to study the Gardner equation from the point of view of Lie symmetries, exact solutions and conservation laws. In \cite{JoKh:10} Johnpillai and Kalique obtained the optimal system of one-dimensional subalgebras of the Lie symmetry algebras of the class
\begin{equation}\label{kha}u_t + u\,u_x+  B(t)\,u_{xxx} + Q(t)\,u =0,\end{equation}
where the linear damping term $B(t)$ and the dispersion term $Q(t)$ are arbitrary smooth functions of the time variable $t$. Later, the same authors \cite{JoKh:11} constructed conservation laws for equation (\ref{kha}) for some special forms of the functions $B(t)$ and $Q(t)$. Vaneeva et al. \cite{Sopho:14} considered the variable coefficient Gardner equations 
\begin{equation}\label{vaneeva}u_t + A(t)\,u\,u_x+ C(t)\,u^{2}u_x + B(t)\,u_{xxx}=0,\end{equation}
where $A(t)$, $B(t)$ and $C(t)$ are smooth functions verifying $B \cdot C \neq 0$. This equation was previously studied by Molati and Ramollo \cite{MoRa:12} who obtained Lie symmetries of equation (\ref{vaneeva}). Vaneeva {\it et al.} enhanced the classification of Lie symmetries obtained in \cite{MoRa:12} through the use of the general extended equivalence group. This enables them to clasify exhaustively these equations by reducing class (\ref{vaneeva}) to the subclass
$$u_t + A(t)\,u\,u_x+ C(t)\,u^{2}u_x + u_{xxx}=0.$$

In this paper, we broaden out the previous results by considering the generalized variable-coefficient Gardner equation with nonlinear terms of any order
\begin{equation}\label{ed1}u_t + A(t)\, u^n\,u_x+ C(t)\,u^{2n}u_x + B(t)\,u_{xxx} + Q(t)\,u =0,\end{equation}
where $n$ is a positive constant, $A(t)$, $B(t) \neq  0$, $C(t) \neq  0$ and $Q(t) \neq 0$
are arbitrary smooth functions of $t$.\\

Due to the variable coefficients that equation (\ref{ed1}) involves, one expects that there exist a transformation by which equation (\ref{ed1}) could be mapped to another equation with the same differential structure but with constant coefficients. An equivalence transformation is a such one.\\

An equivalence transformation is a non-degenerate transformation acting on dependent and independent variables which maps equation (\ref{ed1}) to another equation of the same family, except maybe the form of the coefficients. The main advantage of equivalence transformations is that instead of considering individual equations they permit an analysis for complete equivalent classes. Thus, equivalence transformations appear as a powerful method to study PDEs with variable coefficients.\\

The symmetry group of a PDE is the largest transformation group that acts on dependent and independent variables of the equation so it transforms solutions of the equation in other solutions. Lie symmetry groups is considered to be one of the most powerful methods when analysing PDE. Symmetry groups have several well-known applications. For instance, they can be used to obtain exact solutions \cite{avdonina,Shirvani,wang} or construct conservation laws \cite{Bozhkov,Bru:14,RGB:15,FRE:13,FRE:14,Tra:14}.\\

Given a PDE, a conservation law is a space-time divergence expression which vanishes on all solutions of the PDE. This concept has its origin in physics but their applications spread into many other areas of science. In mathematics, they can be used in numerical methods and mathematical analysis, particularly, investigation of existence, uniqueness and stability of solutions of PDEs. Furthermore, the existence of a large number of conservation laws of a PDE is a strong indicator of its integrability.\\

This work is organised as follows. In Section \ref{equivalence transformations}, we obtain the continuous equivalence group of equation (\ref{ed1}).  Next, in Section \ref{symmetries} we obtain Lie symmetries of the reduced equation obtained by using equivalence transformations. In Section \ref{adjoint equation} and \ref{nonlinearly} we use the concept of adjoint equation and we determine the subclasses of the equation which are nonlinearly self-adjoint. In Section \ref{conservation}, we obtain conservation laws via a general theorem proved by Ibragimov \cite{Ibra:07} and a direct method proposed by Anco and Bluman \cite{anco1, anco}. The conclusions are presented in Section \ref{conclusions}.

\section{Equivalence transformations}\label{equivalence transformations}

In this section we determine the equivalence transformation of class (\ref{ed1}). These transformations allow us to reduce class (\ref{ed1}) to a subclass with simpler form, for instance, reducing the number of arbitrary elements. An equivalence transformation of class (\ref{ed1}) is a nondegenerate point transformation, $\left(t,x,u \right)$ to $\left(\tilde{t},\tilde{x},\tilde{u} \right)$ in the augmented space $\left( t,x,u,A,B,C,Q,n  \right)$ which transforms any equation of class (\ref{ed1}) to an equation from the same class but with different arbitrary elements, $\tilde{A}(\tilde{t})$, $\tilde{B}(\tilde{t})$, $\tilde{C}(\tilde{t})$, $\tilde{Q}(\tilde{t})$ and $\tilde{n}$ from the original ones. We apply Lie's infinitesimal criterion \cite{ovsian} to obtain equivalence transformation of class (\ref{ed1}), i.e., we require not only the invariance of class (\ref{ed1}) but also the invariance of the auxiliary system
\begin{equation}\label{aux}
A_x = A_u=B_x=B_u=C_x=C_u=Q_x=Q_u=n_t=n_x=n_u=0.
\end{equation}

\noindent We consider the one-parameter group of equivalence transformations in $\left( t,x,u,A,B,C,Q,n \right)$ given by
\begin{equation}\label{trans}\begin{array}{rcl} \tilde{t} & = & t+\epsilon \, \tau(t,x,u)+O(\epsilon^2),
\\ \tilde{x} & = & x+\epsilon \, \xi(t,x,u)+O(\epsilon ^2),\\ \tilde{u} & = & u+\epsilon \,
\eta(t,x,u)+O(\epsilon ^2),

\\ \tilde{A} & = & A+\epsilon \,
\omega^1( t,x,u,A,B,C,Q,n )+O(\epsilon ^2),
\\ \tilde{B} & = & B+\epsilon \,
\omega^2( t,x,u,A,B,C,Q,n )+O(\epsilon ^2),
\\ \tilde{C} & = & C+\epsilon \,
\omega^3( t,x,u,A,B,C,Q,n )+O(\epsilon ^2),
\\ \tilde{Q} & = & Q+\epsilon \,
\omega^4( t,x,u,A,B,C,Q,n )+O(\epsilon ^2),
\\ \tilde{n} & = & n+\epsilon \,
\omega^5( t,x,u,A,B,C,Q,n )+O(\epsilon ^2),
\end{array}
\end{equation} where  $\epsilon$ is the
group parameter. In this case, the vector field takes the following form\\
\begin{equation}\label{generator}
 Y  =   \tau \partial_t + \xi \partial_x+\eta \partial_u+ \omega^1 \partial_A  + \omega^2 \partial_B + \omega^3 \partial_C + \omega^4 \partial_Q + \omega^5 \partial_n.   \end{equation}

\noindent Invariance of system (\ref{ed1})-(\ref{aux}) under the one-parameter group of equivalence transformation (\ref{trans}) with infinitesimal generator (\ref{generator}) leads to a system of determining equations. After having solved the determining system, omitting tedious calculations, we obtain the associated equivalence algebra of class (\ref{ed1}) which is infinite dimensional and is spanned by

$$\begin{array}{rcl}\nonumber Y_1 & = & x \partial_x+ u \partial_u + (1-n) A \partial_A + 3 B \partial_B + (1-2n) C \partial_C, \\ \\

\nonumber Y_2 & = & \partial_x,\\ \\

\nonumber Y_\alpha & = & \alpha \partial_t- \alpha_t A \partial_A- \alpha_t B \partial_B- \alpha_t C \partial_C - \alpha_t Q \partial_Q,\\ \\

\nonumber Y_r & = &-r u \partial_u + n r A \partial_A + 2 n r C \partial_C+ r_t \partial_Q,

\end{array}$$

\noindent where $\alpha=\alpha(t)$, $r=r(t)$ are arbitrary smooth functions with $\alpha_t \neq 0$. 

\begin{theorem}\label{th1} The equivalence group of class (\ref{ed1}) consists of the transformations\\
$$ \displaystyle \tilde{t}  =  \alpha(t), \quad \tilde{x}  =  (x+\epsilon_2)e^{\epsilon_1}, \quad \tilde{u}  =  e^{\epsilon_1 - \epsilon_r r(t)}u,$$

$$ \tilde{A}  =  \displaystyle \frac{   e^{n \epsilon_r r+ (1-n)\epsilon_1}}{\alpha_t}A,  \quad   \tilde{B}  =  \displaystyle \frac{e^{3 \epsilon_1}}{\alpha_t}B, \quad
 \tilde{C}   =   \displaystyle \frac{e^{2 n \epsilon_r r+ (1-2n)\epsilon_1}}{\alpha_t}C,  \quad  \tilde{Q}  = \displaystyle \frac{Q+ \epsilon_r r_t }{\alpha_t} \quad\tilde{n}  = n,$$
 
\noindent where $\epsilon_1$, $\epsilon_2$ and $\epsilon_r$ are arbitrary constants. 
\end{theorem}
 
\noindent From Theorem \ref{th1} we obtain that equation (\ref{ed1}) can be transformed into an equation in which the highest order linear term and the nonlinear term of highest order have been set to a nonzero constant values

$$
\tilde{u}_{\tilde{t}}+ \tilde{A}(\tilde{t})\, \tilde{u}^n \, \tilde{u}_{\tilde{x}}+  \tilde{u}^{2 n} \, \tilde{u}_{\tilde{x}}+ \tilde{u}_{\tilde{x}\tilde{x}\tilde{x}}  + \tilde{Q}(\tilde{t}) \, \tilde{u}=0,
$$

\noindent by means of the transformation

\begin{equation}\label{trans1}
\tilde{t}  =  e^{3 \epsilon_1} \int B(t) dt, \quad \tilde{x}  =  (x+\epsilon_2)e^{\epsilon_1}, \quad \tilde{u}  =  \displaystyle e^{\frac{-\epsilon_1}{n}} \left(  \frac{B(t)}{C(t)}  \right)^{-\frac{1}{2 n}} u,
\end{equation}

\noindent where $ \tilde{A}(\tilde{t})= \displaystyle  \frac{e^{-\epsilon_1} A(t)}{\sqrt{B(t) \, C(t)}}$, $ \tilde{B}(\tilde{t})= 1$, $ \tilde{C}(\tilde{t})= 1$ and $ \displaystyle \tilde{Q}(\tilde{t})= e^{-3 \epsilon_1} \left( \frac{Q(t)}{B(t)} + \frac{C(t)}{2 n B(t)^2} \left( \frac{B(t)}{C(t)} \right)_t   \right)  $. Similarly, equation (\ref{ed1})  might be reduced to the same form, for instance, with $\tilde{A}(\tilde{t})=1$. This allow us to restrict without loss of generality our study to the class
\begin{equation}\label{eqtrans}
u_t + A(t)\,u^n \, u_x+ u^{2n} \, u_x + u_{xxx} + Q(t) \, u =0,
\end{equation}

\noindent  due to symmetries and conservation laws obtained for class (\ref{eqtrans}) can be prolonged to class (\ref{ed1}) by using (\ref{trans1}).


\section{Classical Symmetries of class (\ref{eqtrans})}
\label{symmetries}

To apply the Lie classical method to equation (\ref{eqtrans}) we consider the one-parameter Lie group of infinitesimal
transformations in $(t,x,u)$ given by
\begin{eqnarray}\label{ee1}   t^*&=&t+\epsilon \, \tau(t,x,u)+O(\epsilon ^2),
\\\label{ee2} x^*&=&x+\epsilon \, \xi(t,x,u)+O(\epsilon^2),\\\label{ee3} u^*&=&u+\epsilon
\, \eta(t,x,u)+O(\epsilon ^2),\end{eqnarray}
where  $\epsilon$ is the
group parameter. We require that this transformation leaves
invariant the set of solutions of equation (\ref{eqtrans}). This yields
to an overdetermined, linear system of equations for the
infinitesimals $\tau(t,x,u)$, $\xi(t,x,u)$ and $\eta(t,x,u).$ A general element of the symmetry algebra of (\ref{eqtrans}) has the form
\begin{equation}\label{vect1}{\bf v}=
\tau(t,x,u)\partial_t+ \xi(t,x,u) \partial_x+\eta(t,x,u)\partial_u.\end{equation}  Having determined the infinitesimals, the symmetry variables are found by solving the characteristic equation which is equivalent to solving the invariant surface condition
\begin{eqnarray} \nonumber \eta(t,x,u)-\tau(t,x,u)u_t-\xi(t,x,u)u_x=0.\end{eqnarray}

\noindent The set of solutions of equation (\ref{ed1}) is invariant
under the transformation (\ref{ee1})-(\ref{ee3}) provided that
$$pr^{(3)}{\bf v}(\Delta)=0\quad \mbox{when}\quad \Delta=0,$$
where $ \Delta= u_t + A(t)u^n u_x+ u^{2n}u_x + u_{xxx} + Q(t) u$, and $pr^{(3)}{\bf v}$ is the third prolongation of the
vector field (\ref{vect1}) which is given by
\begin{equation}\label{prolong}  pr^{(3)}{\bf v}  ={\bf v} + \zeta^t \partial u_t+ \zeta^x \partial u_x +\zeta^{xxx} \partial u_{xxx}, \end{equation} 
\noindent where $$\zeta^J(t,x,u^{(3)})=D_J(\eta-\tau u_t-\xi u_x)+\tau
u_{Jt} +\xi u_{Jx},$$
with $J=(j_1,\ldots,j_k)$,  $1\leq j_k\leq 2$, $1\leq k\leq 3$, and $u^{(3)}$ denote the sets of partial derivatives up to third order \cite{olver}.
Applying (\ref{prolong}) we obtain a set of determining equations for the infinitesimals. Simplifying this system we get, $\tau=\tau(t)$, $\xi=\xi(x,t)$ and $\eta=\eta(x,t,u)$, where $\tau$, $\xi$ and $\eta$ must satisfy the following conditions:

\begin{equation}\label{sis}\begin{array}{r}

\eta_{uu}=0,\\

\eta_{ux} -\xi_{xx}=0,\\



\tau_t -3 \xi_x =0,\\

\tau u Q_t- \eta_u u Q+ 3 \xi_x u Q +\eta Q
+ \eta_x u^n A +\eta_x u^{2n} + \eta_{xxx} + \eta_t =0, \\

\tau u^{n+1} A_t + 2 \xi_x u^{n+1} A + n \eta u^n A+ 2 \xi_x u^{2n+1} + 2n \eta u^{2n} +3 \eta_{uxx} u - \xi_{xxx} u -\xi_t u =0.\\

\end{array} \end{equation}

\noindent In order to find Lie symmetries of the equation (\ref{eqtrans}) is neccesary to distinguish according to powers of $u$ which appears in (\ref{sis}). From determining system (\ref{sis}), if $n$, $A(t)$ and $Q(t)$ are arbitrary, we obtain ${\bf v}_1= \partial_x$. For the following cases, we obtain new symmetries:\\

\subsection{{\bf Case 1}: $n \neq  1$, $A(t)= \displaystyle a \,(b \, t +c)^{-\frac{1}{3}}$, $Q(t)=d \, (b \, t+c)^{-1}$}
\medskip

\noindent For this special form of $A(t)$ and $Q(t)$ we get the following infinitesimal generators

\begin{eqnarray} \nonumber
{\bf v}_1, \qquad {\bf v}_2= \left( b \, t + c \right) \partial_t+\frac{b}{3}x \partial_x - \frac{b}{3 n}u \partial_u.
\end{eqnarray}

\subsection{{\bf Case 2}: $n=1$, $A(t)= \displaystyle a \,(b \, t +c)^{-d}$, $Q(t)=b \,d \, (b \, t+c)^{-1}$}
\medskip
\noindent In this case, we obtain

\begin{eqnarray} \nonumber
{\bf v}_1, \qquad {\bf v}_2^\prime= \left( b \, t + c \right) \partial_t+ \left( \frac{b}{3}x +\beta \right) \partial_x + \left(\gamma- \frac{b}{3} \right) \partial_u,
\end{eqnarray}

\noindent where $\gamma(t)$ and $\beta(t)$ are given by

\begin{eqnarray}{}
\label{gammac2} \gamma(t)&=&\frac{\left( 3 d-1  \right) a \, b  }{6 }\left(b \,t+c\right)^{-d}, \\ \nonumber\\
\label{betac2} \beta(t) &=&  \left\{ \begin{array}{lcl}
\displaystyle \frac{\left( 3 d-1  \right) a^2  }{6 (1-2d)}\left(b \, t+c\right)^{1-2d}  & \mbox{ si } & \displaystyle d \neq \frac{1}{2}, \\
& & \\
\displaystyle \frac{a^2 \log \left( b \, t+c \right)}{12} & \mbox{ si } & \displaystyle d=\frac{1}{2}.
\end{array}
\right.
\end{eqnarray}


\subsection{\bf Case 3: $n=\displaystyle{\frac{1}{2}}$, $A(t)=0$, $Q(t) \mbox{ \textnormal{arbitrary}}$}
\medskip
\noindent Now, we get the following generators

\begin{eqnarray} \nonumber
{\bf v}_{\tau}= \tau \, \partial_t+  \frac{\tau_t}{3}x \, \partial_x + \left(\frac{\tau_{tt}}{3}x- \frac{2\tau_t}{3}u \right)\, \partial_u, \qquad {\bf v}_{\beta}= \beta \, \partial_x + \beta_t \partial_u,
\end{eqnarray}

\noindent where $\tau=\tau(t)$ and $\beta = \beta(t)$ are given by

\begin{eqnarray}{}
\label{tauc3} \tau(t) &=& \displaystyle  \frac{a}{\sqrt{2 Q(t)^2+Q_t(t)}},\\
\label{betac3} \beta(t) &=& \displaystyle  b \int e^{\int -Q(t) dt} dt+ c.
\end{eqnarray}

\noindent We can observe that in this case, we have obtained two infinitesimal generators with infinite dimension, and therefore we obtain infinite nontrivial conservation laws.\\

\noindent In the above cases, $a$, $b$, $c$ and $d \neq 0$ are arbitrary constants.

\section{Formal Lagrangian and adjoint equation}
\label{adjoint equation}

In order to construct conservation laws we need to set some prior concepts. We consider an sth-order partial differential equation
 \begin{equation}
 \label{fa}
 F({\rm x},u,u_{(1)}, \ldots,u_{(s)}) =0,
 \end{equation}
and 
\begin{equation}
 \label{lag}
 {\cal L}=v\, F({\rm x},u,u_{(1)}, \ldots,u_{(s)}),
 \end{equation}
its formal Lagrangian, where $v=v({\rm x})$ is a new dependent variable, ${\rm x}=(x^1,\ldots,x^n)$ denotes the set of independent variables, $u$ the dependent variable and 
$u_{(1)}=\{u_i\},$ $u_{(2)}=\{u_{ij}\},\ldots\,$ the sets of the
partial derivatives  of the first, second, etc. orders,
$u_i=\partial u/\partial x^i$, $u_{ij}=\partial^2 u/\partial
x^i\partial x^j.$ Let 

$$ \frac{\delta }{\delta u}=\frac{\partial}{\partial u}
 +\displaystyle\sum_{s=1}^{\infty} (-1)^s D_{i_1}\cdots
D_{i_s}\frac{\partial}{\partial u_{i_1\cdots i_s}},$$

\noindent be the variational derivative (the Euler-Lagrange operator). Here
  $$
 D_i=\displaystyle\frac{\partial}{\partial
x^i}+u_i\frac{\partial}{\partial u}+ u_{ij}\frac{\partial}{\partial
u_j}+\cdots
 $$
are the total differentiations. The adjoint equation to {\rm(\ref{fa})} is
 \begin{equation}
 \label{faadj}
 F^{*}({\rm x},u,v,u_{(1)}, v_{(1)},\ldots,u_{(s)},v_{(s)}):=\frac{\delta(v\,F)}{\delta u}
 =0. \end{equation}

\begin{theorem}\label{th2}
The adjoint equation to equation (\ref{eqtrans}) is
\begin{equation}\label{adjointTR} F^*\equiv v Q-u^{2n} v_{x} -v_{x x x} -u^{n} v_{x} A-v_{t}.\end{equation}
\end{theorem}

\section{Nonlinearly self-adjoint equations}
\label{nonlinearly}
We use the following definition given in \cite{Ibran:11}. 

\begin{definition}
Equation {\rm (\ref{fa})} is said to be {\bf nonlinearly self-adjoint} if the
 equation obtained from the adjoint equation {\rm (\ref{faadj})}
 by the substitution \begin{equation}\label{newv}v=\varphi({\rm x},u),\end{equation}
 with a certain function $\varphi({\rm x},u)\neq 0$
 is identical with the original equation {\rm(\ref{fa})}, i.e.,
 \begin{equation}\label{condadj}
 F^*\,{|}_{v=\varphi}=\lambda ({\rm x},u,...) F,\end{equation}
\noindent for some differential function $ \lambda= \lambda({\rm x},u,...)$. If $\varphi=u$  or $\varphi=\varphi(u)$ and $\varphi'(u)\neq 0$, equation {\rm (\ref{fa})} is said {\bf self-adjoint} or {\bf quasi self-adjoint}, respectively. Furthermore, if $\varphi=\varphi({\rm x},u)$ such us $\varphi_u \neq 0$ and $\varphi_{\rm x} \neq 0$, equation (\ref{fa}) is said to be {\bf weak self-adjoint} {\normalfont \cite{weak}}. \end{definition}

\noindent Taking into account the expression (\ref{adjointTR}) and using
(\ref{newv}) and its derivatives we rewrite equation (\ref{condadj}) as

\begin{equation}\label{eqlambda}\begin{array}{rl}
\left(\lambda+\varphi_{u}\right) u_{xxx}+ 3 \varphi_{u u} u_x u_{xx} +3 \varphi_{u x} u_{xx} + \left(   u^n A \lambda + u^{2n} \lambda + \varphi_u u^n A + \varphi_u u^{2n} +3 \varphi_{uxx} \right) u_x &\\ \\
+3 \varphi_{uux} u_x^2 + \varphi_{uuu} u_x^3 +\left(\lambda+\varphi_{u}\right) u_t +u Q \lambda-\varphi Q+\varphi_{x} u^n A+ \varphi_{x} u^{2n}+\varphi_{xxx} +\varphi_{t} & = 0.
\end{array} \end{equation} 
 
\noindent Equation (\ref{eqlambda}) should be satisfied identically in all variables $u_t$, $u_{x}$, $u_{xx},\ldots$. Equating to zero the
coefficients of the derivatives of $u$ we obtain:

\begin{theorem}\label{thnonlinearly}
 Equation (\ref{eqtrans}) with any functions $A(t)$ and $Q(t)$ is nonlinearly self-adjoint and $\varphi(t,x,u)$ is given by
\begin{eqnarray} \nonumber \varphi=p(t)u+q(t,x),\end{eqnarray}
in the following cases:
\begin{itemize}
\item[$\bullet$] If $ \displaystyle n=\frac{1}{2}$ and $A=0$, we have
$$\displaystyle p(t)= c_1 e^{2 \int Q dt}- c_2 e^{2 \int Q dt} \int e^{- \int Q dt} dt , \quad  q(t,x)=\left( c_2 x+ c_3  \right)e^{\int Q dt}. $$

\item[$\bullet$] In any other case, we obtain
$$\displaystyle p(t)=c_1 e^{2\int Q dt}, \quad  q(t)=c_2 e^{\int Q dt}. $$
\end{itemize}
With $c_1$, $c_2$ and $c_3$ are arbitrary constants.
\end{theorem}


\section{Conservation Laws}
\label{conservation}

Conservation laws appear in many of physical, chemical and
mechanical processes, such laws enable us to solve problems in which
certain physical properties do not change in the course of time
within an isolated physical system. In mathematics, the importance of conservation laws lies in that they are strongly related with the integrability of a partial differential equation and can be used to obtain exact solutions, among others.\\ 

\noindent A conservation law of equation (\ref{eqtrans}) is a space-time divergence such that
\medskip
\begin{equation}\label{cl}
 \begin{array}{rcc}
\displaystyle{D_t T(t,x,u,u_x,u_t,...)+D_x X(t,x,u,u_x,u_t,...)=0,}\\
\end{array}
\\
\end{equation}\\
on all solutions $u(t,x)$ of equation (\ref{eqtrans}). Here, $T$ represents the conserved density and $X$ the associated flux \cite{ancoibra}, and $D_x$, $D_t$ denote the total derivative operators with respect to $x$ and $t$ respectively.\\

In this section, we obtain conservation laws of equation (\ref{eqtrans}) using two different methods, the general theorem on conservation laws proved by Ibragimov and the direct method of the multipliers by Anco and Bluman.

\subsection{Conservation laws by using the general theorem by Ibragimov}  
\medskip

In \cite{Ibra:07} Ibragimov showed a new way for obtaining conserved quantities associated with any symmetry of a given differential system wherein the number of equations is equal to the number of dependent variables. To apply this theorem, the concepts of adjoint equation and nonlinear self-adjointness previously defined in Section \ref{adjoint equation} and \ref{nonlinearly} are necessary. The conserved vector obtained by using the theorem proved in \cite{Ibra:07} provides a nonlocal conservation law for the equation. According to this theorem, conserved vectors of equation (\ref{eqtrans}) can be obtained using the following:\\

\begin{theorem}\label{thCL} Any Lie point, Lie-B\"{a}cklund or non-local symmetry
\begin{eqnarray} \nonumber
 {\bf v}=\tau(t,x,u,u_{(1)},\ldots)\frac{\partial}{\partial t}+\xi(t,x,u,u_{(1)},\ldots)\frac{\partial}{\partial x}
 +\eta(t,x,u,u_{(1)},\ldots)\frac{\partial}{\partial u},
 \end{eqnarray} of  equation {\rm (\ref{eqtrans})}
 provides a conservation law $\,D_t T+ D_x X=0\, $ for the simultaneous system {\rm (\ref{eqtrans})}, {\rm
 (\ref{adjointTR})}. The conserved vector is given by
  \begin{eqnarray} \nonumber
\begin{array}{ll}
T&=\tau{\cal L}+ W\left[\displaystyle\frac{\partial{\cal L}}{\partial u_t}\right],\\[3.5ex]
 
 X &= \xi{\cal L}+ W\left[\displaystyle \frac{\partial{\cal L}}{\partial u_x}+D_x^2 \left(\frac{\partial {\cal L}}{\partial u_{xxx}}\right)\right] 
 - D_x \left( W \right) \displaystyle D_x\left( \frac{\partial {\cal L}}{\partial u_{xxx}} \right)   
+D_x^2\left(W\right)  \displaystyle \frac{\partial {\cal L}}{\partial u_{xxx}},
 \end{array}
 \end{eqnarray}
 where ${\cal L}$ is given by (\ref{lag}) and $W$ is defined as follows:
\begin{eqnarray} \nonumber
 W = \eta - \tau u_t-\xi u_x.
 \end{eqnarray}
 \end{theorem}

\noindent Now, we construct conservation laws for the different cases considered in Section \ref{symmetries}.\\

\noindent $\bullet$ From Case 1, we consider generator 

\begin{eqnarray} \nonumber
{\bf v}_2+k \, {\bf v}_1= \left(b t + c \right)  \partial_t + \left( \displaystyle \frac{b}{3}x + k \right)  \partial_x - \frac{b}{3 n}u \, \partial_u,
\end{eqnarray}
\smallskip

\noindent where $k$ is an arbitrary constant. From Theorem \ref{thnonlinearly}, we have
\begin{eqnarray} \nonumber \varphi= \displaystyle{ c_1 \displaystyle \left( b \, t +c \right)^{ \frac{2 d}{b}} u+ c_2 \left( b \, t +c \right)^{ \frac{d}{b}}}.
\end{eqnarray}
\medskip

\noindent Thus, we obtain the conservation law (\ref{cl}) with the conserved vector

$$\begin{array}{lll}T &=& \displaystyle \left(b\,t+c\right)^{{{d}\over{b}}}\,u \left( {{c_{1}\,\left(6\,d\,n+b\,n-2\,b\right)\over{6\,n}} \,\left(b\,t
 +c\right)^{{{d}\over{b}}}\,u}+ {{c_{2}\,\left(3\,d\,n+b\,n-b\right)}\over{3\,n}}  \right),\end{array}$$

$$\begin{array}{lll} X&=& \displaystyle \frac{\left(b\,t
 +c\right)^{{{d}\over{b}}}}{3\,n} \left( {{c_{1}\,\left(6\,d\,n+b\,n-2\,b\right)\,\left(2\,u\,u_{xx}-
 u_{x}^2\right)}\over{2}} \,\left(b\,t
 +c\right)^{{{d}\over{b}}} + {{c_{2}\,\left(3\,d\,n+b\,n-b\right)\,u_{xx}}} \right. \\ \\

& & \displaystyle \left. +{{c_{1}\,\left(6\,d\,n+b\,n-2\,b\right)\over{
 2\,n+2}} \,\left(b\,t
 +c\right)^{{{d}\over{b}}}\,u^{2\,n+2}} + {{c_{2}\,\left(3\,d\,n+b\,n-b\right)}\over{2
 \,n+1}}\,u^{2\,n+1} \right.  \\ \\

& & \displaystyle  \left. +{{c_{1}\,a\,\left(6\,d\,n+b\,n-2\,b\right)}\over{n+2}}\,\left(b\,t+c\right)^{{{d}\over{b}}-{{1}\over{3}}}\,u^{n+
 2} + {{c_{2}\,a\,\left(3\,d\,n+b\,n-b\right)
 }\over{n+1}} \,\left(b\,t+c\right)^{-{{1}\over{3}}}\,u^{n+1} \right)

 .\end{array}$$\\

\noindent $\bullet$ From Case 2, we consider generator

\begin{eqnarray} \nonumber
{\bf v}_2^\prime= \left(b \, t + c \right)  \partial_t + \left( \displaystyle \frac{b}{3}x + \beta \right)  \partial_x + \left( \gamma - \frac{b}{3}u \right) \, \partial_u,
\end{eqnarray}
\smallskip

\noindent where $\gamma$ and $\beta$ are given by (\ref{gammac2}) and (\ref{betac2}) respectively. From Theorem \ref{thnonlinearly}, we have $\varphi$ is given by 
\begin{eqnarray} \nonumber \varphi= c_1 \displaystyle \left( b \, t +c \right)^{ 2 d} u+ c_2 \left( b \, t +c \right)^{ d}.
\end{eqnarray}

\noindent Thus, we obtain the conservation law (\ref{cl}) with the conserved density

$$\begin{array}{lll}T &=& \displaystyle \frac{b \, \left(b\,t+c\right)^{d} \, u }{6} \left( {{c_{1}\,\left(6\,d-1\right)}} \left(b\,t+c\right)^{
 d} u +  {{  c_{1}\,a \left( 3\,d-1\right)+6\,c_{2}\,d}}   \right) +{{c_{2}\,a\,b\,\left(3\,d-1\right)}\over{6}}.\end{array}$$
 
\noindent In the case that $\displaystyle d \neq \frac{1}{2}$ we get

$$\begin{array}{lll} X&=& \displaystyle \frac{b \, \left(b\,t+c\right)^{d}}{6} \left(   c_{1}\,\left(6\,d-1\right)\,\left(2\,u\,u_{xx}-u_{x}^2 + \frac{u^4}{2}\right) \left(b\,t+c\right)
 ^{d}  + \left(  c_{1}\,a\left( 3\,d-1 \right)+6\,c_{2}\,d
 \right)\,u_{xx}   \right. \\ \\

& & \displaystyle \left.   +  \left( c_{1}\,{ a}\left( 5\,d-1 \right)+2\,c_{2}\,
 d\right)\,u^3 \right) + {{a\,b\,\left( c_{1}\,a\left( 3\,d-1 \right)+6\,c_{2}\,d \right)\,u^2}\over{12}} .\end{array}$$
 
\noindent If $\displaystyle d = \frac{1}{2}$ the flux is given by

$$\begin{array}{lll} X&=& \displaystyle \frac{b \, \sqrt{b\,t+c}}{12} \left(   2 \,  c_{1}\,\left(4\,u\,u_{xx}-2\,u_{x}^2 + u^4\right)  \sqrt{b\,t+c}   + \left( c_{1}\,a +6\,c_{2} \right)\,u_{xx} + \left(3\, c_{1}\,a +2\,c_{2} \right)\,u^{3} \right)+ \frac{a \, b\, \left(  c_{1}\,a +6\,c_{2} \right) \, u^2}{24} .\end{array}$$\\
\medskip
\noindent $\bullet$ From Case 3, we consider generator 

\begin{eqnarray} \nonumber
{\bf v}_{\tau}+{\bf v}_{\beta}= \tau  \partial_t + \left( \displaystyle \frac{\tau_t}{3}x + \beta \right)  \partial_x + \left( \frac{\tau_{tt}}{3}x - \frac{2 \tau_t}{3}u+ \beta_t \right) \, \partial_u,
\end{eqnarray}
\smallskip

\noindent where $\tau(t)$ and $\beta(t)$ are given by (\ref{tauc3}) and (\ref{betac3}) respectively. From Theorem \ref{thnonlinearly}, we obtain
\begin{eqnarray} \nonumber \varphi= e^{2 H} \left( c_1-c_2 L \right) u+ e^H (c_2 x +c_3),
\end{eqnarray}

\noindent  where 

\begin{equation}\label{HL} H(t)=\int Q(t) dt, \qquad L(t)=\displaystyle \int  e^{-\int Q(t) dt} dt, \end{equation} 

\noindent with $Q(t)$ arbitrary function.\\

\noindent Thus, we get the conservation law (\ref{cl}) with the conserved vector

$$\begin{array}{lll}T &=& \displaystyle e^{2 H} \left( c_1-c_2 L \right) u \left( \left( \tau Q -\frac{\tau_t}{2} \right) u+   \frac{\tau_{tt}x}{3}+ \beta_t \right)+ e^H (c_2 x +c_3) \left( \tau Q u + \frac{\tau_{tt}x}{3}+ \beta_t \right) -e^H u \left( \frac{c_2 \tau u}{2}- c_2 \, \beta + \frac{c_3 \tau_t}{3} \right),\end{array}$$

$$\begin{array}{lll} X&=& \displaystyle e^{2 H} \left( c_1-c_2 L \right) \left( \tau  Q  \left( 2 u u_{xx}- u_x^2+ \frac{2 u^3}{3}  \right) +  u_{xx} \left( \beta_t -\tau_t u \right)+ \frac{\tau_{tt}x}{6}\left( 2 u_{xx}+u^2 \right)  +u_x \left( \frac{\tau_t}{2}u_x- \frac{\tau_{tt}}{3}  \right)  \right. \\ \\

& & \displaystyle \left. -u^2 \left( \frac{\tau_t}{3}u-\frac{\beta_t}{2} \right) \right)  + \tau e^H Q (c_2 x +c_3)\left(u_{xx}+\frac{u^2}{2}\right)+ \frac{e^H ( 2 u u_{xx}+u^2)}{6}\left( c_1 \tau_{tt}e^H x +3 c_2 \, \beta-c_3 \tau_t \right) \\ \\

& & \displaystyle -c_2 \tau e^H (Q u_x +u u_{xx})+\frac{c_2 e^H}{6}\left( 3 \tau u_x^2-2 \tau u^3-2 \tau_{tt} \right). \end{array}$$\\


\subsection{Conservation laws by using the direct method of the multipliers of Bluman and Anco}
\medskip

Anco and Bluman \cite{anco1, anco} proved a general direct construction method to obtain conservation laws for PDEs which can be expressed in a standard Cauchy-Kovaleskaya form
$$u_t=G(x,u,u_x,u_{xx},\ldots,u_{nx}).$$

\noindent  This method is based on the concept of multiplier. We call multiplier to a certain function $\Lambda(t,x,u,u_x,u_{xx},...)$ which holds that $(u_t+A \, u^n \, u_x+ u^{2n}\, u_x+u_{xxx}+ Q \, u)\Lambda$ is a divergence expression for all functions $u(t,x)$, not just solutions of the equation (\ref{eqtrans}). We suppose, without loss of generality, that $T$, $X$ and the multiplier $\Lambda$ have no dependence on $u_t$ and all derivatives of $u_t$.\\ 

\noindent Using equation (\ref{eqtrans}), we construct an equivalent conservation law in which has been removed $u_t$ and its derivatives from the conserved vector
\begin{eqnarray*} 
 \widehat{T}= T \mid_{u_t=\Delta}= T- \Phi,  \\ \\
 \widehat{X}= X \mid_{u_t=\Delta}= X- \Psi,  
\end{eqnarray*}

\noindent where $ \Delta= - A \, u^n \, u_x-  u^{2n} \, u_x- u_{xxx} -Q \, u$, so that
$$ \left(D_t \widehat{T}(t,x,u,u_x,u_{xx},...)+D_x \widehat{X}(t,x,u,u_x,u_{xx},...)\right)\mid_{u_t=\Delta}=0, $$
and where
\begin{eqnarray*}{} D_t \mid_{u_t=\Delta} &=& \partial_t+ \Delta\partial_u +D_{x}(\Delta)\partial_{u_x}+...    \\ \\
 D_x \mid_{u_t=\Delta} & = &  \partial_x+ u_x \partial_u+ u_{xx} \partial_{u_x}+...=D_x,\\ 
 \end{eqnarray*}

\noindent is held on all solutions of the equation (\ref{eqtrans}). In particular, moving off of solutions, we have the identity

\begin{eqnarray*} D_t = D_t \mid_{u_t=\Delta} + (u_t + A \, u^n \, u_x+ u^{2n}\,u_x + u_{xxx} + Q\,u)\partial_u   + D_x(u_t + A\,u^n\,u_x+ u^{2n}\,u_x + u_{xxx} + Q\,u)\partial_{u_x}+... \end{eqnarray*}
\\
Taking into account the expressions given above we construct the characteristic form of the conservation law (\ref{cl})
\begin{equation}\label{char} \begin{array}{r}  D_t \widehat{T}+D_x \left(\widehat{X}+\widehat{\Psi} \right)   = \left( u_t + A \, u^n\, u_x+  u^{2n}\,u_x +  u_{xxx} + Q \,u \right) \Lambda, \end{array}\end{equation}
\normalsize
where
\begin{eqnarray*} \widehat{\Psi}(t,x,u,u_x,u_t,...)= E_{u_x}(\widehat{T})\left( u_t + A\,u^n\,u_x+ u^{2n}\,u_x + u_{xxx} + Q\, u \right) \\ \\ +E_{u_{xx}}(\widehat{T})D_x\left( u_t + A\,u^n\,u_x+  u^{2n}\,u_x + u_{xxx} + Q\, u \right)+... \end{eqnarray*}
is a trivial flux, and where $E_u= \partial_u-D_x \partial_{u_x}+D_x^{2}\partial_{u_{xx}}-...,$ denotes the (spatial) Euler operator with respect to u. By using the characteristic form (\ref{char}) is deduced that each conserved density in the form (\ref{cl}) arise from a multiplier $\Lambda$ of the equation (\ref{eqtrans}) which has no dependence on $u_t$ and its differential consequences. Multipliers $\Lambda$ are obtained by using that the divergence condition must be verify identically

\begin{eqnarray} \nonumber \frac{\delta}{\delta u} \left(  \left(   u_t + A\,u^n\, u_x+  u^{2n}\,u_x + u_{xxx} + Q\,u  \right) \Lambda  \right) =0, \end{eqnarray}\\

\noindent where $ \displaystyle{ \frac{\delta}{\delta u}= \partial_u-D_x \partial_{u_x}-D_t \partial_{u_t}+D_x D_t\partial_{u_{xt}}+D_x^2\partial_{u_{xx}}+...} $, denotes the variational derivative.\\

This yield us to an overdetermined system of equations for $\Lambda$ which is linear in $u_t$, $u_{tx}$, $u_{txx}$,... Equating to zero these coefficients we obtain the equivalent equations

$$ -D_t \Lambda- (A\,u^n+u^{2n})\, D_x \Lambda- \, D_x^3 \Lambda+  Q \Lambda=0,$$
and
\\
$$ \Lambda_u= E_u(\Lambda), \qquad \, \Lambda_{u_x}=-E_u^{(1)}(\Lambda), \qquad \, \Lambda_{u_{xx}}= -E_u^{(2)}(\Lambda), ... $$\\
which are verified for all solutions $u(t,x)$ of equation (\ref{eqtrans}).
\\

\noindent The conserved density is obtained from a multiplier $\Lambda$ by using a standard method \cite{wolf}
$$ T= \int_{0}^{1} d \lambda \, \, u \Lambda (t,x,\lambda u,\lambda u_x,\lambda u_{xx},...).$$
\\
In this case, we consider multipliers up to second order, i.e., $\Lambda(t,x,u,u_x,u_{xx})$. We obtain that equation (\ref{eqtrans}) admits, for $n$, $A(t)$ and $Q(t)$ arbitrary, the following multiplier\\ 

\begin{eqnarray} \nonumber
\Lambda= \tilde{c_1} e^{2 H}u+ \tilde{c_2} e^{H}.
\end{eqnarray}
\medskip

\noindent From this multiplier, we obtain the conserved vector given by
$$ \begin{array}{lll}T &=& \displaystyle \frac{1}{2}\tilde{c_1} e^{2H} u^2  + \tilde{c_2}e^{H} u,\end{array}$$

$$\begin{array}{lll} X&=& \displaystyle \frac{1}{2}\tilde{c_1} e^{2 H} \left( 2 u u_{xx}-u_x^2 \right)+\tilde{c_2} e^H u_{xx}+\tilde{c_1} e^{2 H} \frac{u^{2n+2}}{2n+2}+ \tilde{c_1} e^{2 H} A \frac{u^{n+2}}{n+2}+ \tilde{c_2} e^H \frac{u^{2n+1}}{2n+1} +\tilde{c_2} e^H A \frac{u^{n+1}}{n+1}.\end{array}$$\\

\noindent Furthermore, in the case that $n=\displaystyle \frac{1}{2}$, $A(t)=0$, we get a new multiplier

\begin{eqnarray} \nonumber
\Lambda= \tilde{c_3} \left(  e^{H} x - e^{2 H} L u \right).
\end{eqnarray}
\medskip

\noindent Consequently, we obtain the following conserved vector
$$\begin{array}{lll}T &=& \displaystyle \frac{\tilde{c_3}}{2} \left(  2 e^{H} x u - e^{2 H} L u^2 \right),\end{array}$$

$$\begin{array}{lll} X&=& \displaystyle \frac{\tilde{c_3}}{6} \left( e^{H} x \left(  6 u_{xx} + 3 u^2 \right)- e^{2 H} L \left(2 u^3 - 3 u_x^2+6 u u_{xx} \right) -6 e^{H} u_x \right).\end{array}$$\\
 
\noindent In the above conserved vectors, $\tilde{c_1}$, $\tilde{c_2}$ and $\tilde{c_3}$ are arbitrary constants, $H(t)$ and $L(t)$ are given by (\ref{HL}).


\section{Conclusions}\label{conclusions}

In this work, we performed a study of a generalized variable-coefficient Gardner equation involving many arbitrary smooth functions. The equation considered generalizes substantially many other equations previously study by other authors by means of symmetries and conservation laws \cite{JoKh:10,JoKh:11,MoRa:12,Sopho:14}. We have obtained the equivalence transformation group of equation (\ref{ed1}), which it is an infinite dimensional group. Equivalence group allows us to enhance our study considering a tranformation which leads us to a subclass (\ref{eqtrans}) of equation (\ref{ed1}) with fewer number of arbitrary functions. Lie symmetries of equation (\ref{eqtrans}) have been obtained.\\

Furthermore, we have proved that equation (\ref{eqtrans}) is nonlinearly self-adjoint. In \cite{Zhang}, it was proved that the substitution $\varphi({\rm x},u)$ used in the property of nonlinear self-adjointness is identical to the multiplier $\Lambda({\rm x},u)$. Nontrivial conservation laws have been obtained by using a general theorem given by Ibragimov \cite{Ibra:07} and the direct method of the multipliers of Anco and Bluman \cite{anco}. The conservation laws obtained via Ibragimov's method can be constructed by using Theorem \ref{thCL}, which does not use the integral of functions, in contrast to the multipliers method. It can be seen that although substitutions and multipliers are similar, conservation laws provided from them may be different.

\section*{Acknowledgments}

The authors acknowledge the financial support from Junta de Andaluc\'ia group FQM--201, University of C\'adiz. The second and third author also acknowledge the support of DGICYT project MTM2009-11875 with the participation of FEDER.


\begin{thebibliography}{xx}  

\bibitem{anco1}
S.C.~Anco and G.~Bluman.
\newblock  Direct constrution method for conservation laws of partial differential equations Part I: Examples of conservation law classifications. \newblock  \emph{Euro. Jnl of Applied mathematics}, 13, 545--566, 2002.

\bibitem{anco}
S.C.~Anco and G.~Bluman. 
\newblock  Direct constrution method for conservation laws of partial differential equations Part II: General treatment.
\newblock \emph{Euro. Jnl of Applied mathematics}, 13, 567--585, 2002.

\bibitem{ancoibra}
S.C. Anco et al. 
\newblock   Symmetries and conservation laws of the generalized Krichever-Novikov equation.
\newblock  Preprint, 2014.

\bibitem{avdonina}
E.D.~Avdonina and N.H.~Ibragimov.
\newblock Conservation laws and exact solutions for nonlinear diffusion in anisotropic media.
\newblock \emph{Commun Nonlinear Sci Numer Simulat}, 18, 2595--2603, 2013.

\bibitem{Bozhkov}
Y.~Bozhkov, S.~Dimas and N.H.~Ibragimov.
\newblock Conservation laws for a coupled variable-coefficient modified Korteweg-de Vries system in a two-layer fluid model.
\newblock \emph{Commun Nonlinear Sci Numer Simulat}, 18, 1127--1135, 2013.






\bibitem{Bru:14}
M.S. Bruz\'on, M.L. Gandarias and R. de la Rosa.
\newblock Conservation laws of a  family Reaction-Diffusion-Convection Equations.
\newblock Localized Excitations in Nonlinear Complex Systems, Series: Nonlinear Systems and Complexity, 7, 2014.

\bibitem{RGB:15} R.~de la Rosa, M.~L.~Gandarias and M.~S.~Bruz\'{o}n.
\newblock A study for the microwave heating of some chemical reactions through Lie symmetries and conservation laws.
\newblock \emph{J Math Chem}, 53, 949--957, 2015.

\bibitem{FRE:13}
I.L. Freire.
\newblock  New classes of nonlinearly self-adjoint evolution equations of third-and
fifth-order.
\newblock \emph{Commun Nonlinear Sci Numer Simulat}, 18(3), 493--499, 2013.

\bibitem{FRE:14}
I.L.~Freire and J.C.S.~Sampaio.
\newblock On the nonlinear self-adjointness and local conservation laws for a class of evolution equations unifying many models.
\newblock \emph{Commun Nonlinear Sci Numer Simulatn}, 19, 350--360, 2014.

\bibitem{weak}
M.L.~Gandarias.
\newblock Weak self-adjoint differential equations.
\newblock \emph{J. Phys. A},  44, 262001--262007, 2011.




\bibitem{JoKh:10}
A.G.~Johnpillai and C.M.~Khalique.
\newblock Gruop analysis of KdV equation with time dependent coefficients.
\newblock \emph{Applied Mathematics and Computation}, 216, 3761--3771, 2010.

\bibitem{JoKh:11}
A.G.~Johnpillai and C.M.~Khalique.
\newblock Conservation laws of KdV equation with time dependent coefficients.
\newblock \emph{Commun Nonlinear Sci Numer Simulat}, 16, 3081--3089, 2011.

\bibitem{Ibra:07}
N.H.~Ibragimov.
\newblock A new conservation theorem.
\newblock \emph{J. Math. Anal. Appl.}, 333, 311--328, 2007.

\bibitem{Ibran:11}
N.H.~Ibragimov.
\newblock Nonlinear self-adjointness and conservation laws.
\newblock \emph{J. Phys. A: Math. Theor.}, 44, 432002--432010, 2011.

\bibitem{MoRa:12}
M.~Molati and M.P.~Ramollo.
\newblock Symmetry classification of the Gardner equation with time-dependent coefficients arising in stratified fluids.
\newblock \emph{Commun Nonlinear Sci Numer Simulat}, 17, 1542--1548, 2012.

\bibitem{olver} P.~Olver.
\newblock Applications of Lie groups to differential equations. 
\newblock Springer-Verlag, New York, 1993.
  
\bibitem{ovsian}
L.V.~Ovsyannikov.
\newblock Group analysis of differential equations. 
\newblock Academic, New York, 1982.

\bibitem{Shirvani}  V.~Shirvani-Sh and M.~Nadjafikhah.
\newblock Conservation laws and exact solutions of the Whitham-type equations.
\newblock \emph{Commun Nonlinear Sci Numer Simulat}, 19, 2212--2219, 2014. 



\bibitem{Tra:14}
R.~Tracin\'a, M.S.~Bruz\'on, M.L.~Gandarias and M.~Torrisi.
\newblock  Nonlinear self-adjointness, conservation laws, exact
solutions of a system of dispersive evolution equations.
\newblock \emph{Commun Nonlinear Sci Numer Simulat}, 19, 3036--3043,
2014.

\bibitem{Sopho:14} O.~Vaneeva, O.~Kuriksha and C.~Sophocleous,
\newblock Enhanced group classification of Gardner equations with time-dependent coefficients,
\newblock \emph{Commun Nonlinear Sci Numer Simulat}, 22, 1243--1251, 2015. 

\bibitem{wang} 
G.W.~Wang, X.G.~Liu and Y.Y~Zhang.
\newblock Symmetry reduction, exact solutions and conservation laws of a new fifth-order nonlinear integrable equation.
\newblock \emph{Commun Nonlinear Sci Numer Simulat}, 18, 2313--2320, 2013.

\bibitem{wolf} T.~Wolf,
\newblock An efficieny improved program LIEPDE for determining Lie-symmetries of PDEs.
\newblock {\em Proceedings of Modern Group Analysis: advances analytical and computational methods in mathematical physics}, 377-385, 1993.


\bibitem{Zhang} 
Z.Y.~Zhang,
\newblock On the existence of conservation law multiplier for partial differential equations.
\newblock \emph{Commun Nonlinear Sci Numer Simulat}, 20, 338--351, 2015.




\end{thebibliography}
\end{document}